\documentclass[11pt,reqno]{amsart}
\usepackage{amsmath,amssymb,latexsym, amsthm, amscd, mathrsfs, stmaryrd} 

\usepackage[all]{xy}
\usepackage{hyperref}
\usepackage{color}
\usepackage{bm}
\newcommand{\nc}{\newcommand}
\nc{\browntext}[1]{\textcolor{brown}{#1}}
\nc{\greentext}[1]{\textcolor{green}{#1}}
\nc{\redtext}[1]{\textcolor{red}{#1}}
\nc{\bluetext}[1]{\textcolor{blue}{#1}}
\nc{\brown}[1]{\browntext{ #1}}
\nc{\green}[1]{\greentext{ #1}}
\nc{\red}[1]{\redtext{ #1}}
\nc{\blue}[1]{\bluetext{ #1}}

\newcommand{\ev}{\overline{0}} 
\newcommand{\odd}{\overline{1}} 
\newcommand{\Y}{\check{E}}
\newcommand{\kk}{h}
\newcommand{\vs}{\boldsymbol{\varsigma}}
\newcommand{\qvs}{\blue{q\boldsymbol{\varsigma}}}

\newcommand{\dvev}[1]{{{B}_{\ev}^{{(#1)}}}}

\newcommand{\B}{{\mbf U}^\imath}
\newcommand{\dvd}[1]{B_{\odd}^{{(#1)}}}
\newcommand{\dv}[1]{{B}^{{(#1)}}_{\odd}}
\newcommand{\ttt}{B} 

\allowdisplaybreaks

\newtheorem{thm}{Theorem}[section]

\newtheorem{lem}[thm]{Lemma}
\newtheorem{prop}[thm]{Proposition}

\newtheorem{example}[thm]{Example}

\theoremstyle{remark}
\newtheorem{rem}[thm]{Remark}

\numberwithin{equation}{section}

\newcommand{\mbf}{\mathbf}

\newcommand{\N}{\mathbb N}
\newcommand{\Q}{\mathbb Q}
\newcommand{\Z}{\mathbb Z}

\newcommand{\U}{\mbf U}

\newcommand{\Ui}{{\mathbf U}^\imath}

\newcommand{\qbinom}[2]{\begin{bmatrix} #1\\#2 \end{bmatrix} }

\title[Formulae of $\imath$-divided powers III] {Formulae of $\imath$-divided powers in ${\bf U}_q(\mathfrak{sl}_2)$, III
}
\author[Xinhong Chen]{Xinhong Chen}
\address{ Department of Mathematics\\ Southwest Jiaotong University\\Sichuan 611756, China}
\email{chenxinhong@swjtu.edu.cn}
\author[Weiqiang Wang]{Weiqiang Wang}
\address{ Department of Mathematics\\ University of Virginia\\ Charlottesville, VA 22904, USA}
\email{ww9c@virginia.edu}


\begin{document}

	\begin{abstract}
		The $\imath$-divided powers (depending on a parity) form the $\imath$canonical basis for the split rank 1 $\imath$quantum group and they are a basic ingredient for $\imath$quantum groups of higher rank. We obtain closed formulae for the structure constants for multiplication of the $\imath$-divided powers. Closed formulae for the comultiplication of the $\imath$-divided powers are also obtained. These structure constants are integral and positive.
	\end{abstract}
	
	\maketitle
	\setcounter{tocdepth}{1}
	\tableofcontents

	\section{Introduction}

	\subsection{}
	
	The quantum group of $\mathfrak{sl}_2$, $\U = \langle E, F, K^{\pm 1} \rangle$, contains  divided powers $F^{(n)}=F^n/[n]!$. 
	These divided powers form the canonical basis for half a quantum group $\U^- =\Q(q) [F]$, and they are a basic ingredient for higher rank quantum groups as well, cf. \cite{L93}. The divided powers satisfy a simple recursive relation
	\begin{align}
		F\cdot F^{(n)}=[n+1]F^{(n+1)},
	\end{align}
	and admit a closed multiplication formula
	\begin{align}
		\label{MF}
		F^{(m)} \cdot F^{(n)}=\qbinom{m+n}{n} F^{(m+n)}.
	\end{align}
	Under the comultiplication $\Delta$ in \eqref{eq:Delta}, we have
	\begin{align}
		\label{coMF}
		\Delta (F^{(n)}) = \sum_{a=0}^n q^{a(n-a)} F^{(a)} \otimes F^{(n-a)} K^{-a}.
	\end{align}

	\subsection{}
	
	The $\imath$quantum groups arise from the theory of quantum symmetric pairs (see \cite{Let99}). Fix $\vs\in \Q(q)^\times$.
	The $\imath$quantum group $\Ui$ (of split rank 1)  is a coideal $\Q(q)$-subalgebra of $\U$ \cite{K93}, which is the polynomial algebra in  $\ttt$:
	\[
	\Ui =\Q(q)[\ttt],
	\]
	where
	\begin{equation}  \label{eq:t}
		\ttt:= F + \vs EK^{-1}.
	\end{equation}
	The {\em $\imath$-divided powers} for $\vs=q^{-1}$ are $\imath$canonical basis for the $\imath$quantum group $\Ui$ \cite{BW18a}. Explicit formulae for $\imath$-divided powers for $\vs=q^{-1}$ were conjectured therein (for a somewhat different $B$; see Remark~\ref{rem:B}), and the conjecture was established in \cite{BeW18}. The same formulae of $\imath$-divided powers were shown  in \cite{BeW21} to be valid, for more general $\ttt= F + q^{-1} EK^{-1} + [k] K^{-1}$, for $k \in \Z$. (The special cases for $k=0,1$ were treated in \cite{BeW18}, and the case for $k=1$ appeared first in \cite{BW18a}.)
	
	The $\imath$-divided powers for $\vs=q^{-1}$ have played a fundamental role in the general theory of $\imath$canonical bases \cite{BW21}.  The $\imath$-divided powers with a general parameter $\vs$ appears in the Serre presentation of $\imath$quantum groups \cite{CLW21} and also in Serre-Lusztig relations and braid group action for $\imath$quantum groups \cite{CLW21b}.
	
	The (split) $\imath$-divided powers come in 2 forms $\{\dvev{n} \mid n \ge 0 \}$ and $\{\dv{n} \mid n \ge 0 \}$, depending on a parity $\ev, \odd$ (of the highest weight of a finite dimensional simple $\U$-module). The formulae for the $\imath$-divided powers $\dvev{n}$ (respectively, $\dv{n}$) can be found in \eqref{def:idpev} (respectively, \eqref{def:idpodd}).
	
	The $\imath$-divided powers $\dvev{n}$  are determined by $\dvev{0}=1$ and the following recursive relations,  for $a\ge 1$,
	\begin{align}
		\label{eq:ttev}
		\begin{split}
			\ttt \cdot \dvev{2a-1} &=[2a] \dvev{2a},
			\\
			\ttt \cdot \dvev{2a} &=  [2a+1] \dvev{2a+1} +\qvs [2a] \dvev{2a-1}.
		\end{split}
	\end{align}
	For example,
	$
	\dvev{1} =\ttt,
	\dvev{2} =\ttt^2/[2],$ and
	$\dvev{3} = \ttt(\ttt^2-\qvs[2]^2)/[3]!.$
	
	There are similar recursive formulae for the $\imath$-divided powers $\dv{n}$; see \eqref{eq:ttodd}.

	\subsection{}
	
	The goal of this paper is to determine closed formulae for multiplication and comultiplication of $\imath$-divided powers. These formulae are more involved than their counterparts \eqref{MF}--\eqref{coMF} for Lusztig's divided powers, since $\imath$-divided powers are not monomials in $B$. It is remarkable that such closed formulae actually exist. The harder part of this work is to find these closed formulae, especially the comultiplication formulae. The proofs are routine by sometimes lengthy induction.

	$\triangleright$ The multiplication formulae for the $\imath$-divided powers $\dvev{n}$ are given in Theorem~\ref{thm:evmulti}.
	
	$\triangleright$ The multiplication formulae for $\dv{n}$ are given in Theorem~\ref{thm:oddmulti}.
	
	$\triangleright$ The comultiplication formulae for $\dvev{n}$ are given in Theorem~\ref{thm:even}.
	
	$\triangleright$ The comultiplication formulae for $\dv{n}$ are given in Theorem~\ref{thm:odd}.
	
	The comultiplication formulae are reminiscent of the PBW expansion formulae for the $\imath$-divided powers; cf. \cite{BeW18, BeW21}.

	\subsection{}
	Assume $\vs=q^{-1}$.
	The structure constants for multiplication of $\imath$-divided powers (see Theorem~\ref{thm:evmulti} and Theorem~\ref{thm:oddmulti}) are typically given by ratios of products of quantum integers. It may not be straightforward (though it is an interesting combinatorial problem) to recognize that they are integral and positive, i.e., they lie in $\N[q,q^{-1}]$. As $\imath$-divided powers are integral basis for an $\imath$quantum group over $\Z[q,q^{-1}]$ \cite{BW18a}, these structure constants indeed lie in $\Z[q,q^{-1}]$. Positivity of these structure constants follows from the geometric interpretation of $\imath$canonical basis of type AIII \cite{LW18, Li20}. Similarly, the structure constants for comultiplication with respect to ($\imath$-)canonical bases are positive \cite{FL21}; in the rank 1 setting, it follows from our closed formulae (see Remark~\ref{rem:positive}).
	
	New $\imath$-divided powers associated to various rank 1 $\imath$quantum groups were constructed in \cite{BW21} with favorable integrality property, but otherwise remain to be poorly understood (in contrast to the ones studied in this paper); it is desirable to study and understand them better.
	
	The formulae obtained in this work will be applied in some ongoing work on braid group action. It shall also find applications to $\imath$quantum groups at roots of unity (as we learned from Huanchen Bao). These are the main reasons for which we decide to make the notes public (the main results have been known to us for some time). 
	We hope they may be helpful in $\imath$categorification as well.
	
	%

\vspace{3mm}
{\bf Acknowledgement.}
We thank Collin Berman for his help in 2018 with Mathematica computation regarding comultiplication formulae.
XC is supported by the  National Natural Science Foundation of China (No. 12271447), and the Fundamental Research Funds for the Central
Universities grant 2682020ZT100 and 2682021ZTPY043.
WW is partially supported by the NSF grant DMS-2001351.

	\section{Multiplication formulae for $\dvev{n}$}
	
	In this section, we recall the $\imath$-divided powers, and then compute the multiplication formulae for $\imath$-divided powers $\dvev{n}$.

	\subsection{$\imath$-Divided powers}
	
	Let $\U$ be the quantum group of $\mathfrak{sl}_2$ over $\Q(q)$ with standard generators $E,F, K^{\pm1}$.
	There is a comultiplication $\Delta: \U \rightarrow \U \otimes \U$ \cite{L93} such that
	\begin{equation}
		\label{eq:Delta}
		\Delta(E)  = E \otimes 1 +  K \otimes E, \quad \Delta(F) = 1 \otimes F + F \otimes K^{-1}, \quad \Delta(K) = K  \otimes K.
	\end{equation}
	We have
	\begin{align}
		\label{D:ttt}
		\Delta(B) =  \ttt \otimes K^{-1} + 1 \otimes (F + \vs EK^{-1}).
	\end{align}
	
	Denote by $\N$ the set of nonnegative integers. 
	The following $\imath$-divided powers $\dvev{n}$ and $\dv{n}$ appear in \cite{CLW21}, which generalize the ones in \cite{BW18a, BeW18} for $\vs =q^{-1}$.  
	
	For $n\in \N$, let 
	\begin{align}
		B_{\ev}^{(n)} &= \frac{1}{[n]!}\left\{ \begin{array}{ccccc} B\prod_{j=1}^k (B^2-\qvs[2j]^2 ) & \text{if }n=2k+1,\\
		\prod_{j=1}^{k} (B^2-\qvs[2j-2]^2) &\text{if }n=2k. \end{array}\right.
\label{def:idpev}
\end{align}
See  also \eqref{eq:ttev} for a recursive definition.
			%

	
	For $n\in \N$, let 
\begin{align}
	B_{\odd}^{(n)} &=\frac{1}{[n]!}\left\{ \begin{array}{ccccc} B\prod_{j=1}^k (B^2-\qvs [2j-1]^2 ) & \text{if }n=2k+1,\\
		\prod_{j=1}^k (B^2-\qvs [2j-1]^2) &\text{if }n=2k; \end{array}\right.
	\label{def:idpodd}
\end{align}
	
	The $\imath$-divided powers $\dv{n}$ are determined by the following recursive relations: $\dv{0}=1$, and
	\begin{align}  
		\label{eq:ttodd}
		\begin{split}
			\ttt \cdot \dv{2a} &=  [2a+1] \dv{2a+1} ,
			\\
			\ttt  \cdot \dv{2a+1} &=[2a+2] \dv{2a+2}+\qvs [2a+1] \dv{2a}.
		\end{split}
	\end{align}
	
	\begin{rem}
		\label{rem:B}
		The above formulae for $\imath$-divided powers (with $\boldsymbol{\varsigma}=q^{-1}$) remain unchanged if we consider more general $\ttt= F + q^{-1} EK^{-1} +[2p] K^{-1}$ (for fixed $p\in \Z$) other than \eqref{eq:t}  with $\boldsymbol{\varsigma}=q^{-1}$; see \cite{BeW21}.
		
		The $\imath$-divided powers of $\dot B :=F+ q^{-1} EK^{-1} + [2p+1]K^{-1}$ (for fixed $p\in \Z$) are given by the same formulae above (with $\boldsymbol{\varsigma}=q^{-1}$) but with the parities $\ev$ and $\odd$ swapped; see \cite{BeW21}; such formulae first appeared in \cite{BW18a}, which correspond to $\dot B$ with $p=0$.
	\end{rem}

	\subsection{The multiplication formulae}
	
	The following simple identities can be verified directly.
	\begin{lem}
		For $n, m, \ell \in \Z$, we have
		\begin{align}
			[n+m]  +[n-m] &= [n] [2]_{q^m}, \label{eq:eq1}
			\\
			[n+m] [n-m] &= [n]^2 -[m]^2, \label{eq:eq2}
			\\
			[m] [m+n] -[\ell][\ell+n] &= [m-\ell] [m+\ell+n], \label{eq:eq3}
			\\
			[2n] &= [2] [n]_{q^2} \label{eq:eq4}.
		\end{align}
	\end{lem}

	We have the following multiplication formulae for $\imath$-divided powers $\dvev{n}$.
	
	\begin{thm}
		\label{thm:evmulti}
		For $k, a \ge 1$, we have
		\begin{align}  
			& \dvev{2k-1}   \dvev{2a-1} =
			\label{oddodd} \\
			& \qbinom{2k +2a -2}{2k-1}
			\left( \dvev{2k+2a-2}+ \sum_{\ell=2}^k \prod_{m=2}^\ell \frac{  [2a-2m+2] [2k-2m+2]}{ [2k+2a-2m+1]  [2m-2]} (\qvs)^{\ell-1} \dvev{2k+2a-2\ell} \right),
			\notag
			\\
			& \dvev{2k-1}  \dvev{2a} =
			\label{oddev} \\
			&\qbinom{2k +2a -1}{2k-1}
			\left( \dvev{2k+2a-1}+ \sum_{\ell=1}^k \prod_{m=1}^\ell \frac{  [2a-2m+2] [2k-2m+2]}{ [2k+2a-2m+1]  [2m]}(\qvs)^\ell \dvev{2k+2a-2\ell-1} \right),
			\notag
			\\  
			& \dvev{2k}  \dvev{2a-1} =
			\label{evodd} \\
			&\qbinom{2k +2a -1}{2k}
			\left( \dvev{2k+2a-1}+ \sum_{\ell=1}^k \prod_{m=1}^\ell \frac{  [2a-2m+2] [2k-2m+2]}{ [2k+2a-2m+1]  [2m]} (\qvs)^\ell\dvev{2k+2a-2\ell-1} \right),
			\notag
			\\
			&\dvev{2k} \dvev{2a} =
			\label{evev} \\
			& %
			\qbinom{2k +2a}{2k}
			\left( \dvev{2k+2a}+ \sum_{\ell=1}^k  \frac{[2k+2a-2\ell]}{[2k+2a]} \prod_{m=1}^\ell \frac{  [2a-2m+2] [2k-2m+2]}{ [2k+2a-2m+1]  [2m]}(\qvs)^\ell \dvev{2k+2a-2\ell} \right).
			\notag
		\end{align}
	\end{thm}

	\begin{example}
		
		\begin{align*}
			\dvev{2}   \dvev{2a-1} &=  \qbinom{2a+1}{2} \dvev{2a+1} +\frac{[2a]^2}{[2]} (\qvs)\dvev{2a-1},
			\\
			\dvev{2}  \dvev{2a} &=  \qbinom{2a+2}{2} \dvev{2a+2} +\frac{[2a]^2}{[2]} (\qvs)\dvev{2a},
			\\
			\dvev{3}   \dvev{2a-1} &=  \qbinom{2a+2}{3} \dvev{2a+2} + \frac{[2a+2][2a][2a-2]}{[3]!} (\qvs) \dvev{2a},
			\\
			\dvev{3}  \dvev{2a} &=  \qbinom{2a+3}{3} \dvev{2a+3}  + \frac{[4]}{[2]}\qbinom{2a+2}{3}  (\qvs) \dvev{2a+1}  +  \frac{[2a+2][2a][2a-2]}{[3]!} (\qvs)^2\dvev{2a-1},
			\\
			\dvev{4}  \dvev{2a-1} &=  \qbinom{2a+3}{4} \dvev{2a+3}
			+  \frac{[2a+2][2a+1][2a]^2}{[3][2]^2} (\qvs) \dvev{2a+1}
			+  \frac{[2a+2][2a]^2[2a-2]}{[4]!}(\qvs)^2 \dvev{2a-1},
			\\ %
			\dvev{4}  \dvev{2a} &=  
			\qbinom{2a+4}{4}  \dvev{2a+4}
			+   \frac{[2a+2]^2[2a+1][2a]}{[3][2]^2} (\qvs)  \dvev{2a+2}
			+  \frac{[2a+2][2a]^2[2a-2]}{[4]!}(\qvs)^2 \dvev{2a}.
		\end{align*}
		
	\end{example}

	\subsection{Proof of Theorem \ref{thm:evmulti}}
	
	Note  in each of the four formulae \eqref{oddodd}--\eqref{evev}
	above the summation up to $k$ can be replaced by the (same) summation up to $\min (k,a)$, since the additional terms clearly vanish when $k>a$. Therefore if we switch $a$ and $k$ then
	Equations~\eqref{oddodd} and \eqref{evev} remain unchanged while \eqref{oddev} and \eqref{evodd} get swapped. We have chosen to keep both equivalent formulations in the theorem to facilitate the inductive proof below.
	
	
	We prove by induction on $k$. Below let us add indices to mark the identities in the theorem as \eqref{oddodd}$_{k}$, \eqref{oddev}$_k$, \eqref{evodd}$_k$, \eqref{evev}$_k$. For the base cases of induction, the identities \eqref{evodd}$_0$--\eqref{evev}$_0$ are trivial, while the identities \eqref{oddodd}$_1$--\eqref{oddev}$_1$ are simply \eqref{eq:ttev}.
	
	The proof of the theorem will be completed in the following steps (i)-(iv):
	\begin{enumerate}
		\item[(i)]
		\eqref{oddev} $\Leftrightarrow$ \eqref{evodd};
		\item[(ii)]
		\eqref{oddev}$_{k}$ $\Rightarrow$ \eqref{evev}$_{k}$;
		\item[(iii)]
		\eqref{evodd}$_{k}$+\eqref{oddodd}$_{k}$   $\Rightarrow$ \eqref{oddodd}$_{k+1}$;
		\item[(iv)]
		\eqref{oddodd}$_{k}$ $\Rightarrow$ \eqref{evodd}$_{k}$.
	\end{enumerate}
	Step (i) follows by swapping $k$ and $a$ as we already noted, while (ii) follows by applying \eqref{eq:ttev} directly.
	
	Let us prove (iii): \eqref{evodd}$_{k}$+\eqref{oddodd}$_{k}$   $\Rightarrow$ \eqref{oddodd}$_{k+1}$.
	By \eqref{eq:ttev}, we have
	\begin{align*}
		&\dvev{2k+1}   \dvev{2a-1} =\frac{\ttt\cdot \dvev{2k}-(\qvs)[2k]\dvev{2k-1}}{[2k+1]} \dvev{2a-1}\\
		=&\frac{\ttt}{[2k+1]} \qbinom{2k +2a -1}{2k}
		\Big( \dvev{2k+2a-1}+ \sum_{\ell=1}^k \prod_{m=1}^\ell \frac{  [2a-2m+2] [2k-2m+2]}{ [2k+2a-2m+1]  [2m]} (\qvs)^\ell\dvev{2k+2a-2\ell-1} \Big)\\
		&-\frac{[2k](\qvs)}{[2k+1]}  \qbinom{2k +2a -2}{2k-1}
		\Big( \dvev{2k+2a-2}+ \sum_{\ell=2}^k \prod_{m=2}^\ell \frac{  [2a-2m+2] [2k-2m+2]}{ [2k+2a-2m+1]  [2m-2]} (\qvs)^{\ell-1}\dvev{2k+2a-2\ell} \Big)\\
		=& \qbinom{2k +2a}{2k+1}
		\left( \dvev{2k+2a}+ \sum_{\ell=1}^k \prod_{m=1}^\ell \frac{[2k+2a-2\ell]  [2a-2m+2] [2k-2m+2]}{ [2k+2a-2m+1][2k+2a]  [2m]} (\qvs)^\ell\dvev{2k+2a-2\ell} \right.\\
		&-\frac{[2k]^2(\qvs)}{[2k+2a][2k+2a-1]}  
		\left.\Big( \dvev{2k+2a-2}+ \sum_{\ell=2}^k \prod_{m=2}^\ell \frac{  [2a-2m+2] [2k-2m+2]}{ [2k+2a-2m+1]  [2m-2]} (\qvs)^{\ell-1}\dvev{2k+2a-2\ell} \Big)\right)
		\\
		%
		=& \qbinom{2k +2a}{2k+1}\left(
		\dvev{2k+2a}+ \sum_{\ell=2}^{k+1} \prod_{m=2}^{\ell} \frac{[2k+2a-2\ell+2]  [2a-2m+4] [2k-2m+4]}{ [2k+2a-2m+3][2k+2a]  [2m-2]}(\qvs)^{\ell-1}\right.
		\\
		&\quad\quad\cdot \dvev{2k+2a-2\ell+2}-\frac{[2k]^2}{[2k+2a][2k+2a-1]}
		\\
		& \quad\quad\cdot
		\left.\Big( (\qvs)\dvev{2k+2a-2}+ \sum_{\ell=3}^{k+1} \prod_{m=2}^{\ell-1} \frac{  [2a-2m+2] [2k-2m+2]}{ [2k+2a-2m+1]  [2m-2]} (\qvs)^{\ell-1}\dvev{2k+2a-2\ell+2} \Big)\right)\\
		=& \qbinom{2k +2a}{2k+1}\Big(
		\dvev{2k+2a}+ \sum_{\ell=2}^{k+1} \prod_{m=2}^{\ell} \frac{[2a-2m+2] [2k-2m+4]}{[2a-2\ell+2] [2k+2a-2m+3][2k+2a]  [2m-2]}\\
		&\quad \quad \quad \quad \quad  \cdot\big([2k+2a-2l+2][2a]-[2k][2l-2]\big)
		(\qvs)^{\ell-1} \dvev{2k+2a-2\ell+2} \Big).
	\end{align*}
	The last formula can be easily rewritten as \eqref{oddodd}$_{k+1}$ by noting that
	\begin{align*}
		[2k+2a-2l+2][2a]-[2k][2l-2]=[2a-2\ell+2][2k+2a].
	\end{align*}

	Now we prove (iv):  \eqref{oddodd}$_{k}$ $\Rightarrow$ \eqref{evodd}$_{k}$.
	By \eqref{eq:ttev}, we have
	\begin{align*}
		\dvev{2k} & \dvev{2a-1}=\frac{\ttt\cdot \dvev{2k-1}}{[2k]}\dvev{2a-1}\\
		=&\frac{\ttt}{[2k]}\qbinom{2k +2a -2}{2k-1}
		\Big( \dvev{2k+2a-2}+ \sum_{\ell=2}^k \prod_{m=2}^\ell \frac{  [2a-2m+2] [2k-2m+2]}{ [2k+2a-2m+1]  [2m-2]} (\qvs)^{\ell-1}\dvev{2k+2a-2\ell} \Big)\\
		=&\qbinom{2k +2a -1}{2k}
		\Big(\dvev{2k+2a-1}+ \frac{[2k+2a-2]}{[2k+2a-1]}(\qvs)\dvev{2k+2a-3}\\
		&  +\sum_{\ell=2}^k \prod_{m=2}^\ell \frac{[2k+2a-2\ell+1]  [2a-2m+2] [2k-2m+2]}{[2k+2a-1] [2k+2a-2m+1]  [2m-2]} (\qvs)^{\ell-1}\dvev{2k+2a-2\ell+1} \\
		&  +\sum_{\ell=2}^k \prod_{m=2}^\ell \frac{[2k+2a-2\ell]  [2a-2m+2] [2k-2m+2]}{[2k+2a-1] [2k+2a-2m+1]  [2m-2]} (\qvs)^\ell\dvev{2k+2a-2\ell-1}
		\Big)
		\\
		%
		=&\qbinom{2k +2a -1}{2k}
		\Big(\dvev{2k+2a-1}+ \frac{[2k+2a-2]}{[2k+2a-1]}(\qvs)\dvev{2k+2a-3}\\
		& +\sum_{\ell=1}^{k-1} \prod_{m=1}^{\ell} \frac{[2k+2a-2\ell-1]  [2a-2m] [2k-2m]}{[2k+2a-1] [2k+2a-2m-1]  [2m]} (\qvs)^\ell\dvev{2k+2a-2\ell-1} \\
		& +\sum_{\ell=2}^k \prod_{m=2}^\ell \frac{[2k+2a-2\ell]  [2a-2m+2] [2k-2m+2]}{[2k+2a-1] [2k+2a-2m+1]  [2m-2]} (\qvs)^\ell \dvev{2k+2a-2\ell-1}
		\Big) \\
		%
		=&\qbinom{2k +2a -1}{2k}
		\Big(\dvev{2k+2a-1}+ \frac{[2][2k+2a-2]+[2a-2][2k-2]}{[2k+2a-1][2]}(\qvs)\dvev{2k+2a-3}\\
		&  +\sum_{\ell=2}^{k-1} \prod_{m=1}^{\ell} \frac{[2k+2a-2\ell-1]  [2a-2m+2] [2k-2m+2][2a-2\ell][2k-2\ell]}{[2k+2a-2m+1]  [2m]
			[2a] [2k][2k+2a-2\ell-1]} (\qvs)^\ell\dvev{2k+2a-2\ell-1} \\
		&  +\sum_{\ell=2}^k \prod_{m=1}^\ell \frac{[2k+2a-2\ell]  [2a-2m+2] [2k-2m+2][2\ell]}{ [2k+2a-2m+1]  [2m][2a][2k]} (\qvs)^\ell\dvev{2k+2a-2\ell-1}
		\Big)
		\\
		%
		=&\qbinom{2k +2a -1}{2k}
		\Big(\dvev{2k+2a-1}+\sum_{\ell=1}^{k-1} \prod_{m=1}^{\ell} \frac{[2a-2m+2] [2k-2m+2]}{[2k+2a-2m+1] [2m]
			[2a] [2k]} \\
		&\big([2a-2\ell][2k-2\ell]+ [2\ell][2k+2a-2\ell] \big)(\qvs)^\ell\dvev{2k+2a-2\ell-1} \\
		&+ \prod_{m=1}^k \frac{  [2a-2m+2] [2k-2m+2]}{ [2k+2a-2m+1]  [2m]}(\qvs)^k \dvev{2a-1}
		\Big).
	\end{align*}
	The last formula can be easily rewritten as \eqref{evodd}$_{k}$ by 
	applying the identity \eqref{eq:eq3}. 
	This completes the proof of Steps (i)--(iv) and hence of Theorem~\ref{thm:evmulti}.

	\section{Multiplication formulae for $\dvd{n}$}

	In this section, we present the multiplication formulae for $\imath$-divided powers $\dvd{n}$. 

	\subsection{The multiplication formulae for $\dvd{n}$}
	
	\begin{thm}
		\label{thm:oddmulti}
		For $k,a\ge 0$, we have
		\begin{align}  
			& \dv{2k} \dv{2a}  =
			\label{evev2}
			\textstyle
			\qbinom{2k +2a}{2k}
			\left( \dv{2k+2a}+ \sum\limits_{\ell=1}^k  
			\prod\limits_{m=1}^\ell \frac{  [2a-2m+2] [2k-2m+2]}{ [2k+2a-2m+1]  [2m]} (\qvs)^\ell\dv{2k+2a-2\ell} \right),
			\\
			&\dv{2k}  \dv{2a+1}  =
			\label{evodd2}
			\textstyle \qbinom{2k +2a +1}{2k}
			\left( \dv{2k+2a+1}+ \sum\limits_{\ell=1}^k \prod\limits_{m=1}^\ell \frac{  [2a-2m+2] [2k-2m+2]}{ [2k+2a-2m+3]  [2m]} (\qvs)^\ell\dv{2k+2a-2\ell+1} \right),
			\\
			&\dv{2k+1}  \dv{2a}= \label{oddev2}
			\textstyle \qbinom{2k +2a +1}{2k+1}
			\left( \dv{2k+2a+1}+ \sum\limits_{\ell=1}^{k }\prod\limits_{m=1}^\ell \frac{  [2a-2m+2] [2k-2m+2]}{ [2k+2a-2m+3]  [2m]} (\qvs)^\ell\dv{2k+2a-2\ell+1} \right),
			\\
			&\dv{2k+1}  \dv{2a+1}=
			\label{oddodd2}
			\textstyle
			\qbinom{2k +2a+2 }{2k+1}
			\left( \dv{2k+2a+2}
			+
			\sum\limits_{\ell=1}^{k+1} \Big(\frac{[2k +2a-2 \ell+2 ][2k-2\ell +2]}{[2k +2a+2 ][2k+2]}\right.
			\\
			&\textstyle
			\qquad\qquad\qquad\quad
			\left.\quad +\frac{[2k +2a-2\ell+3]^2[2\ell]}{[2k +2a+2 ][2a-2\ell+2 ][2k+2]}\Big)\cdot \prod\limits_{m=1}^\ell \frac{  [2a-2m+2] [2k-2m+4]}{ [2k+2a-2m+3]  [2m]} (\qvs)^\ell\dv{2k+2a-2\ell+2}
			\right).
			\notag
		\end{align}
	\end{thm}
	
	\begin{example}
		\begin{align*}
			\dv{2}   \dv{2a} &=  \qbinom{2a+2}{2} \dv{2a+2} +(\qvs)\frac{[2a+2] [2a]}{[2]} \dv{2a},
			\\
			\dv{2}   \dv{2a+1} &=  \qbinom{2a+3}{2} \dv{2a+3} +\frac{[2a+2] [2a]}{[2]}(\qvs) \dv{2a+1},
			\\
			\dv{3}  \dv{2a} &=  
			\qbinom{2a+3}{3} \dv{2a+3}  + (\qvs)\qbinom{2a+2}{3} \dv{2a+1},
			\\
			\dv{3}   \dv{2a+1} &
			=  \qbinom{2a+4}{3} \dv{2a+4} + \frac{[2a+2]^2[2a] +  [2a+3]^2[2a+2]}{[3]!} (\qvs)\dv{2a+2}
			\\
			&\qquad\qquad+ \qbinom{2a+2}{3} (\qvs)^2\dv{2a},
			%
			\\
			\dv{4}  \dv{2a} &=  \qbinom{2a+4}{4} \dv{2a+4} + \frac{[2a+4]}{[2]} \qbinom{2a+2}{3} (\qvs)\dv{2a+2}\\
			&\qquad\qquad+\frac{[2a+4][2a+2][2a][2a-2]}{[4]!}(\qvs)^2\dv{2a},
			\\
			%
			\dv{4}   \dv{2a+1} 
			&= \qbinom{2a+5}{4}  \dv{2a+5} + \frac{[2a]}{[2]} \qbinom{2a+4}{3}(\qvs) \dv{2a+3}
			\\
			&\qquad\qquad+\frac{[2a+4][2a+2][2a][2a-2]}{[4]!}(\qvs)^2\dv{2a+1}.
		\end{align*}
	\end{example}

	\subsection{Proof of Theorem \ref{thm:oddmulti}}
	
	We prove by induction on $k$. Below let us add indices to mark the identities of Theorem \ref{thm:oddmulti} as \eqref{evev2}$_k$,  \eqref{evodd2}$_k$, \eqref{oddev2}$_k$, \eqref{oddodd2}$_{k}$. Note that \eqref{evev2}$_0$--\eqref{evodd2}$_0$ are trivial while \eqref{oddev2}$_1$, \eqref{oddodd2}$_1$ hold by \eqref{eq:ttodd}.
	
	The proof of the theorem will be completed in the following steps.
	\begin{enumerate}
		\item[(i)]
		\eqref{evodd2}  $\Leftrightarrow$ \eqref{oddev2};
		\item[(ii)]
		\eqref{evodd2}$_{k}$ $\Rightarrow$ \eqref{oddodd2}$_{k}$;
		\item[(iii)]
		\eqref{oddodd2}$_{k}+$ \eqref{evodd2}$_{k}$ $\Rightarrow$ \eqref{evodd2}$_{k+1}$;
		\item[(iv)]
		\eqref{oddev2}$_{k}+$ \eqref{evev2}$_{k}$ $\Rightarrow$ \eqref{evev2}$_{k+1}$.
	\end{enumerate}
	Step (i) follows by swapping $a$ and $k$.
	
	Let us prove (ii): \eqref{evodd2}$_{k}$ $\Rightarrow$ \eqref{oddodd2}$_{k}$. Using \eqref{eq:ttodd}, we have
	
	\begin{align*}
		&\dv{2k+1}  \dv{2a+1} =\frac{\ttt \cdot \dv{2k} } {[2k+1]}\dv{2a+1}
		\\
		= &\frac {\ttt}{[2k+1]} \qbinom{2k +2a +1}{2k}
		\Big( \dv{2k+2a+1}+ \sum_{\ell=1}^k \prod_{m=1}^\ell \frac{  [2a-2m+2] [2k-2m+2]}{ [2k+2a-2m+3]  [2m]}(\qvs)^\ell \dv{2k+2a-2\ell+1} \Big)\\
		= &\frac {1}{[2k+1]}\qbinom{2k +2a +1}{2k}\Big([2k+2a+2]\dv{2k+2a+2}+[2k+2a+1](\qvs)\dv{2k+2a}\\
		&\quad \quad \quad \quad +  \sum_{\ell=1}^k \prod_{m=1}^\ell \frac{  [2a-2m+2] [2k-2m+2]}{ [2k+2a-2m+3]  [2m]}(\qvs)^\ell\\
		& \quad \quad \quad \quad  \cdot \big([2k+2a-2\ell+2]\dv{2k+2a-2\ell+2}
		+ [2k+2a-2\ell+1](\qvs)\dv{2k+2a-2\ell} \big)\Big)\\
		=&\qbinom{2k +2a +2}{2k+1}\Big(\dv{2k+2a+2}+\frac{[2k+2a+1]}{[2k+2a+2]}(\qvs)\dv{2k+2a}\\
		&\quad \quad \quad \quad+  \sum_{\ell=1}^k \frac{[2k+2a-2\ell+2]}{[2k+2a+2]} \prod_{m=1}^\ell \frac{  [2a-2m+2] [2k-2m+2]}{ [2k+2a-2m+3]  [2m]}(\qvs)^\ell\dv{2k+2a-2\ell+2}\\
		&\quad \quad \quad \quad + \sum_{\ell=1}^k \frac{[2k+2a-2\ell+1]}{[2k+2a+2]} \prod_{m=1}^\ell \frac{  [2a-2m+2] [2k-2m+2]}{ [2k+2a-2m+3]  [2m]}(\qvs)^{\ell+1}\dv{2k+2a-2\ell} \Big)
	\end{align*}
	\begin{align*}
		=&\qbinom{2k +2a +2}{2k+1}\Big(\dv{2k+2a+2}+\frac{[2k+2a+1]}{[2k+2a+2]}(\qvs)\dv{2k+2a}+\frac{[2k+2a]}{[2k+2a+2]} \frac{  [2a] [2k]}{ [2k+2a+1]  [2]}(\qvs)\dv{2k+2a} \\
		&  +  \sum_{\ell=2}^k \frac{[2k+2a-2\ell+2]}{[2k+2a+2]} \prod_{m=1}^\ell \frac{  [2a-2m+2] [2k-2m+2]}{ [2k+2a-2m+3]  [2m]}(\qvs)^{\ell}\dv{2k+2a-2\ell+2}\\
		&  + \sum_{\ell=2}^{k} \frac{[2k+2a-2\ell+3]^2[2\ell]}{[2k+2a+2][2a-2\ell+2][2k-2\ell+2] } \prod_{m=1}^{\ell} \frac{  [2a-2m+2] [2k-2m+2]}{ [2k+2a-2m+3]  [2m]}(\qvs)^\ell\dv{2k+2a-2\ell+2}\\
		& +\frac{[2a+1]}{[2k+2a+2]} \prod_{m=1}^k \frac{  [2a-2m+2] [2k-2m+2]}{ [2k+2a-2m+3]  [2m]}(\qvs)^{k+1}\dv{2a}  \Big)
		\\
		%
		=& \qbinom{2k +2a+2 }{2k+1}
		\left( \dv{2k+2a+2} + \sum_{\ell=1}^{k} \Big(\frac{[2k +2a-2 \ell +2]}{[2k +2a +2]}+\frac{[2k +2a-2\ell+3]^2[2\ell]}{[2k +2a +2][2k-2\ell +2][2a-2\ell +2]} \Big)\right.\\
		&
		\cdot\prod_{m=1}^\ell \frac{  [2a-2m+2] [2k-2m+2]}{ [2k+2a-2m+3]  [2m]} (\qvs)^{\ell}\dv{2k+2a-2\ell+2} \\ \notag
		&
		\left.+ \frac{[2a+1]^2}{[2k+2a+2][2a-2k]}
		\prod_{m=1}^{k+1} \frac{[2a-2m+2] }{ [2k+2a-2m+3]  } (\qvs)^{k+1}\dv{2a}\right)
		\\
		=& \qbinom{2k +2a+2 }{2k+1}
		\left( \dv{2k+2a+2}
		+
		\sum_{\ell=1}^{k+1} \Big(\frac{[2k +2a-2 \ell+2 ][2k-2\ell +2]}{[2k +2a+2 ][2k+2]}\right.
		\\
		&\left.+\frac{[2k +2a-2\ell+3]^2[2\ell]}{[2k +2a+2 ][2a-2\ell+2 ][2k+2]}\Big)\cdot \prod_{m=1}^\ell \frac{  [2a-2m+2] [2k-2m+4]}{ [2k+2a-2m+3]  [2m]} (\qvs)^\ell \dv{2k+2a-2\ell+2}
		\right).
	\end{align*}

	\vspace{2mm}
	
	Let us now prove (iii): \eqref{oddodd2}$_{k}+$ \eqref{evodd2}$_{k}$ $\Rightarrow$ \eqref{evodd2}$_{k+1}$. Using \eqref{eq:ttodd}, we have
	\begin{align*}
		&\dv{2k+2}  \dv{2a+1} =\frac{\ttt \cdot \dv{2k+1}-(\qvs)[2k+1]\dv{2k} } {[2k+2]}\dv{2a+1}\\
		=&\frac{\ttt}{[2k+2]}\dv{2k+1}\dv{2a+1}-\frac{[2k+1]}{[2k+2]}(\qvs)\dv{2k}\dv{2a+1}\\
		=& \qbinom{2k +2a+2 }{2k+1}
		\frac{\ttt}{[2k+2]} \left( \dv{2k+2a+2} +
		\sum_{\ell=1}^{k+1} \Big(\frac{[2k +2a-2 \ell +2][2k-2\ell +2]}{[2k +2a +2][2k+2]}\right.
		\\
		&\left.+\frac{[2k +2a-2\ell+3]^2[2\ell]}{[2k +2a +2][2a-2\ell +2][2k+2]}\Big) \cdot \prod_{m=1}^\ell \frac{  [2a-2m+2] [2k-2m+4]}{ [2k+2a-2m+3]  [2m]} (\qvs)^\ell\dv{2k+2a-2\ell+2}\right)\\
		&-\frac{[2k+1]}{[2k+2]}\qbinom{2k +2a +1}{2k}
		\\
		&\cdot
		\Big( (\qvs)\dv{2k+2a+1}+ \sum_{\ell=1}^k \prod_{m=1}^\ell \frac{  [2a-2m+2] [2k-2m+2]}{ [2k+2a-2m+3]  [2m]} (\qvs)^{\ell+1}\dv{2k+2a-2\ell+1} \Big)
	\end{align*}
	\begin{align*}
		=& \qbinom{2k +2a+2 }{2k+1}
		\frac{1}{[2k+2]} \bigg([2k+2a+3] \dv{2k+2a+3} \\ \notag
		& +
		\sum_{\ell=1}^{k+1} \Big(\frac{[2k +2a-2 \ell +2][2k-2\ell +2]}{[2k +2a +2][2k+2]}+\frac{[2k +2a-2\ell+3]^2[2\ell]}{[2k +2a +2][2a-2\ell +2][2k+2]}\Big)\\
		& \cdot \prod_{m=1}^\ell \frac{  [2a-2m+2] [2k-2m+4]}{ [2k+2a-2m+3]  [2m]}
		[2k+2a-2\ell+3](\qvs)^\ell\dv{2k+2a-2\ell+3}\bigg)\\
		&-\frac{[2k+1]}{[2k+2]}\qbinom{2k +2a +1}{2k}
		\\
		&\cdot\Big( (\qvs)\dv{2k+2a+1}+ \sum_{\ell=1}^k \prod_{m=1}^\ell \frac{  [2a-2m+2] [2k-2m+2]}{ [2k+2a-2m+3]  [2m]} (\qvs)^{\ell+1}\dv{2k+2a-2\ell+1} \Big)
	\end{align*}
	\begin{align*}
		=& \qbinom{2k +2a+3 }{2k+2}
		\bigg( \dv{2k+2a+3}   \\
		& +\sum_{\ell=1}^{k+1} \Big(\frac{[2k +2a-2 \ell +2][2k-2\ell +2]}{[2k +2a +2][2k+2]}+\frac{[2k +2a-2\ell+3]^2[2\ell]}{[2k +2a +2][2a-2\ell +2][2k+2]}\Big)\\
		&  \cdot \frac{[2k+2a-2\ell+3]}{[2k+2a+3]}
		\prod_{m=1}^\ell \frac{  [2a-2m+2] [2k-2m+4]}{ [2k+2a-2m+3]  [2m]}(\qvs)^\ell \dv{2k+2a-2\ell+3}  \\
		& -\frac{[2k+1]^2}{[2k+2a+3][2k+2a+2]}
		\\
		&\cdot\Big((\qvs)\dv{2k+2a+1}+ \sum_{\ell=1}^k \prod_{m=1}^\ell \frac{  [2a-2m+2] [2k-2m+2]}{ [2k+2a-2m+3]  [2m]} (\qvs)^{\ell+1}\dv{2k+2a-2\ell+1}\Big)\bigg)
	\end{align*}
	\begin{align*}
		=& \qbinom{2k +2a+3 }{2k+2}
		\bigg( \dv{2k+2a+3} \\
		& +\sum_{\ell=1}^{k+1} \Big(\frac{[2k +2a-2 \ell +2][2k-2\ell +2]}{[2k +2a +2][2k+2]}+\frac{[2k +2a-2\ell+3]^2[2\ell]}{[2k +2a +2][2a-2\ell +2][2k+2]}\Big)\\
		&\cdot
		\prod_{m=1}^\ell \frac{  [2a-2m+2] [2k-2m+4]}{ [2k+2a-2m+5]  [2m]} (\qvs)^\ell\dv{2k+2a-2\ell+3}
		-\sum_{\ell=1}^{k+1}\frac{[2k+1]^2[2\ell]}{[2k+2a+2][2a-2\ell+2][2k+2]}\\
		&\cdot   \prod_{m=1}^\ell \frac{  [2a-2m+2] [2k-2m+4]}{ [2k+2a-2m+5]  [2m]} (\qvs)^{\ell}\dv{2k+2a-2\ell+3}\bigg)\\
		=& \qbinom{2k +2a+3 }{2k+2}
		\bigg( \dv{2k+2a+3}    \\
		& +\sum_{\ell=1}^{k+1} \Big(\frac{[2k +2a-2 \ell +2][2a-2\ell+2][2k-2\ell +2]+[2k +2a-2\ell+3]^2[2\ell]-[2k+1]^2[2\ell]}{[2k +2a +2][2a-2\ell +2][2k+2]}\Big)\\
		& \cdot
		\prod_{m=1}^\ell \frac{  [2a-2m+2] [2k-2m+4]}{ [2k+2a-2m+5]  [2m]} (\qvs)^\ell\dv{2k+2a-2\ell+3} \bigg).
	\end{align*}
	The RHS above can be converted to RHS \eqref{evodd2}$_{k+1}$ using the following identity (which  can be derived easily via \eqref{eq:eq2}--\eqref{eq:eq3})
	\begin{align}\label{eq:eq5}
		&[2k +2a-2 \ell +2][2a-2\ell+2] [2k-2\ell+2]+ [2k +2a-2\ell+3]^2[2\ell]-[2k+1]^2  [2\ell]\\ \notag
		=&[2a-2\ell+2][2k+2][2k+2a+2 ]. \notag
	\end{align}
	
	\vspace{2mm}
	
	Let us now prove (iv): \eqref{oddev2}$_{k}+$ \eqref{evev2}$_{k}$ $\Rightarrow$ \eqref{evev2}$_{k+1}$. Using \eqref{eq:ttodd}, we have
	
	\begin{align*}
		&\dv{2k+2}\dv{2a}=\frac{\ttt\cdot\dv{2k+1}-[2k+1](\qvs)\dv{2k}}{[2k+2]}\dv{2a}\\
		=&\frac{\ttt}{[2k+2]}\qbinom{2k +2a +1}{2k+1}
		\Big( \dv{2k+2a+1}+ \sum_{\ell=1}^k \prod_{m=1}^\ell \frac{  [2a-2m+2] [2k-2m+2]}{ [2k+2a-2m+3]  [2m]} (\qvs)^\ell\dv{2k+2a-2\ell+1} \Big)\\
		&-\frac{[2k+1]}{[2k+2]}\qbinom{2k +2a}{2k} (\qvs)
		\Big(\dv{2k+2a}+ \sum_{\ell=1}^k
		\prod_{m=1}^\ell \frac{  [2a-2m+2] [2k-2m+2]}{ [2k+2a-2m+1]  [2m]} (\qvs)^\ell\dv{2k+2a-2\ell} \Big)
		\\
		=&\frac{1}{[2k+2]}\qbinom{2k +2a +1}{2k+1}\Big([2k+2a+2]\dv{2k+2a+2}+[2k+2a+1](\qvs)\dv{2k+2a}  \Big)\\
		&+\frac{1}{[2k+2]}\qbinom{2k +2a +1}{2k+1} \sum_{\ell=1}^k [ 2k+2a-2\ell+2 ] \prod_{m=1}^\ell \frac{  [2a-2m+2] [2k-2m+2]}{ [2k+2a-2m+3]  [2m]}(\qvs)^\ell\dv{2k+2a-2\ell+2} \\
		&+\frac{1}{[2k+2]}\qbinom{2k +2a +1}{2k+1} \sum_{\ell=1}^k [2k+2a-2\ell+1]\prod_{m=1}^\ell \frac{  [2a-2m+2] [2k-2m+2]}{ [2k+2a-2m+3]  [2m]} (\qvs)^{\ell+1}\dv{2k+2a-2\ell}\\
		&-\frac{[2k+1]}{[2k+2]}\qbinom{2k +2a}{2k}
		\Big((\qvs) \dv{2k+2a}+ \sum_{\ell=1}^k
		\prod_{m=1}^\ell \frac{  [2a-2m+2] [2k-2m+2]}{ [2k+2a-2m+1]  [2m]} (\qvs)^{\ell+1}\dv{2k+2a-2\ell} \Big)
		\\
		=&\qbinom{2k +2a +2}{2k+2}\bigg(\dv{2k+2a+2}+\frac{[2k+2a][2a][2k]+[2k+2a+1]^2[2]-[2k+1]^2[2]}{[2k+2a+2][2k+2a+1][2] }(\qvs)\dv{2k+2a}\\
		&+\sum_{\ell=2}^k \frac{[ 2k+2a-2\ell+2 ]}{ [ 2k +2a +2]} \prod_{m=1}^\ell \frac{  [2a-2m+2] [2k-2m+2]}{ [2k+2a-2m+3]  [2m]}(\qvs)^\ell\dv{2k+2a-2\ell+2}\\
		&+\sum_{\ell=2}^{k} \frac{[ 2k+2a-2\ell+3 ]^2[2\ell] }{[2a-2\ell+2] [2k-2\ell+2]  [ 2k +2a +2]}\prod_{m=1}^{\ell} \frac{  [2a-2m+2] [2k-2m+2]}{ [2k+2a-2m+3]  [2m]}(\qvs)^\ell\dv{2k+2a-2\ell+2}\\
		&-\sum_{\ell=2}^{k} \frac{[2k+1]^2[2\ell] }{[2a-2\ell+2] [2k-2\ell+2] [ 2k +2a +2]}\prod_{m=1}^\ell \frac{  [2a-2m+2] [2k-2m+2]}{ [2k+2a-2m+3]  [2m]}(\qvs)^\ell\dv{2k+2a-2\ell+2}\\
		&+\frac{[2a+1]^2-[2k+1]^2}{[2k+2a+2][2a+1]} \prod_{m=1}^k \frac{  [2a-2m+2] [2k-2m+2]}{ [2k+2a-2m+3]  [2m]}(\qvs)^{k+1}\dv{2a}\bigg)\\
		=&\qbinom{2k +2a +2}{2k+2}\bigg(\dv{2k+2a+2}\\
		&+\sum_{\ell=1}^k \frac{[ 2k+2a-2\ell+2 ][2a-2\ell+2] [2k-2\ell+2]+[ 2k+2a-2\ell+3 ]^2[2\ell] -[2k+1]^2[2\ell]}{ [ 2k +2a +2][2a-2\ell+2] [2k-2\ell+2]} \\
		&\cdot \prod_{m=1}^\ell \frac{  [2a-2m+2] [2k-2m+2]}{ [2k+2a-2m+3]  [2m]}(\qvs)^\ell\dv{2k+2a-2\ell+2}
		\\
		&+\frac{[2a-2k]}{[2a+1]} \prod_{m=1}^k \frac{  [2a-2m+2] [2k-2m+2]}{ [2k+2a-2m+3]  [2m]}(\qvs)^{k+1}\dv{2a}\bigg).
	\end{align*}
	The RHS above can be easily converted to RHS \eqref{evev2}$_{k+1}$ using  \eqref{eq:eq5}.
	
	This completes the proof of Theorem~ \ref{thm:oddmulti}.

	\section{Comultiplication formulae for $\dvev{n}$}
	
	In this section, we shall establish a closed formula for $\Delta (\dvev{n})$ in two different forms.

	\subsection{An anti-involution}
	
	Set
	\begin{equation*}
		\Y :=\vs EK^{-1},
		\qquad
		\kk :=\frac{K^{-2}-1}{q^2-1}.
	\end{equation*}
	
	Define, for $a\in \Z, n\ge 0$,
	\begin{equation}  \label{kbinom}
		\qbinom{\kk;a}{n} =\prod_{i=1}^n \frac{q^{4a+4i-4} K^{-2} -1}{q^{4i} -1},
		\qquad
		[\kk;a]= \qbinom{\kk;a}{1}.
	\end{equation}
	Then we have, for $a\in \Z, n\in \N$,
	\begin{align}
		\label{FEk}
		F \Y  -q^{-2} \Y  F =\qvs \kk,
		\qquad \qbinom{\kk;a}{n} F = F \qbinom{\kk;a+1}{n},
		\qquad \qbinom{\kk;a}{n} \Y  =\Y  \qbinom{\kk;a-1}{n}.
	\end{align}

	\begin{lem}
		[\text{\cite[Lemma~2.2]{BeW18}}]
		\label{lem:antiinvolution}
		(1) There is an anti-involution $\chi$ on the $\Q$-algebra $\U$ which sends $E \mapsto E, F \mapsto F, K \mapsto K, q \mapsto q^{-1}$.
		
		(2) Assume $\vs=q^{-1}$. Then $\chi$ restricts to an anti-involution of the $\Q$-algebra $\B$ sending
		\begin{equation}
			\label{eq:vs}
			F\mapsto F,\quad q^{-1}EK^{-1} \mapsto q^{-1}EK^{-1}, \quad  K^{-1} \mapsto K^{-1},  \quad  q \mapsto q^{-1}.
		\end{equation}
		
		Moreover, 
		$\chi$ sends
		\begin{align}
			\label{eq:vs-hn}
			\begin{split}
				\kk\mapsto -q^{2} \kk, &
				\quad
				\dvev{n} \mapsto \dvev{n},
				\quad
				\dv{n} \mapsto \dv{n},
				\\
				\qbinom{\kk;a}{n} &\mapsto (-1)^n   q^{2n(n+1)}  \qbinom{\kk;1-a-n}{n}, \; \forall a\in\Z, n\in\N.
			\end{split}
		\end{align}
	\end{lem}

	\subsection{The comultiplication formulae}
	
	For $x\in \mathbb R$, we denote
	\begin{align*}
		\lfloor x\rfloor &=\max\{m\in \mathbb {Z} \mid m\leq x\},
		\qquad \lceil x\rceil =\min\{n\in \mathbb {Z} \mid n\geq x\}.
	\end{align*}
	
	\begin{thm}
		\label{thm:even}
		For $n$ even, we have
		
		\begin{align}  \label{eq:deltaevev}
			\Delta(\dvev{n})  &=
			\sum_{r=0}^n \dvev{n-r}  \otimes   \bigg(\sum_{c=0}^{\lfloor \frac{r}2 \rfloor}\sum_{a=0}^{r-2c} q^{\binom{2c}{2}+(r-2c)(r-n)-a(r-2c-a)} (\qvs)^c\Y^{(a)}   \qbinom{\kk;-\lfloor \frac{r-2}2 \rfloor}{c} K^{r-n} F^{(r-2c-a)} \bigg).
		\end{align}
		For $n$ odd, we have
		\begin{align} \label{eq:deltaevodd}
			\Delta(\dvev{n})  &=
			\sum_{r=0}^n \dvev{n-r}  \otimes   \bigg(\sum_{c=0}^{\lfloor \frac{r}2 \rfloor}\sum_{a=0}^{r-2c}  q^{\binom{2c+1}{2}+(r-2c)(r-n)-a(r-2c-a)} (\qvs)^c \Y^{(a)}  \qbinom{\kk;-\lfloor \frac{r-1}2 \rfloor}{c} K^{r-n} F^{(r-2c-a)} \bigg).
		\end{align}
	\end{thm}

	\begin{example}
		\begin{align*}
			\Delta(\dvev{2}) &= \dvev{2} \otimes K^{-2}  +\ttt \otimes \big(q^{-1} \Y K^{-1} +q^{-1} K^{-1}F\big)
			+1\otimes \big(\Y^{(2)}  + q^{-1} \Y F+ F^{(2)}  +  q (\qvs)[\kk;0] \big),
			\\ %
			\Delta(\dvev{3}) &=  \dvev{3} \otimes K^{-3}  +\dvev{2} \otimes \big(q^{-2} \Y K^{-2} +q^{-2} K^{-2}F\big)
			\\
			&\quad +\ttt \otimes \big(q^{-2} \Y^{(2)} K^{-1} + q^{-3} \Y K^{-1} F+ q^{-2} K^{-1} F^{(2)}  +  q^{3}(\qvs) [\kk;0] K^{-1} \big)
			\\
			&\quad +1 \otimes  \big(\Y^{(3)} + q^{-2} \Y^{(2)}  F+ q^{-2} \Y F^{(2)}  +F^{(3)} +q^{3} (\qvs)\Y [\kk;-1] +q^{3}(\qvs) [\kk;-1]F \big).
		\end{align*}
	\end{example}
	
	\begin{rem}
		\label{rem:positive}
		When passing to the modified quantum group, the terms $\qbinom{\kk; k }{c}$ in \eqref{eq:deltaevev}--\eqref{eq:deltaevodd} are replaced by some quantum integers up to some $q$-powers, and then the formulas have positive integral coefficients. 
	\end{rem}

	\subsection{Proof of Theorem \ref{thm:even}}
	
	Denote by $S_{n,r} \in \U$ (for both $n$ even and odd)  the expressions in parentheses in the theorem such that
	
	\[
	\Delta(\dvev{n}) =\sum_{r=0}^n \dvev{n-r}  \otimes S_{n,r}.
	\]
	We shall prove the theorem by induction on $n$. The base case when $n=1$ is the formula \eqref{eq:Delta}. 
	The induction is carried out in two steps {\bf I--II}.
	
	\vspace{.3cm}
	{\bf Step~I.}
	Assuming the second formula in Theorem \ref{thm:even} for odd $n=2\ell-1$, we shall prove the first formula for even $n=2\ell$.
	It follows from \eqref{eq:ttev} that
	\begin{align}
		\label{D:2l}
		[2\ell] \Delta(\dvev{2\ell}) = \Delta(\ttt) \Delta(\dvev{2\ell-1})
		= \left(\ttt \otimes K^{-1} + 1 \otimes (\Y +  F) \right)\cdot \Delta(\dvev{2\ell-1}).
	\end{align}
	We shall compare the summands $\dvev{2\ell-r} \otimes -$ in both sides of \eqref{D:2l}. As the summands on the RHS of
	\eqref{D:2l} are known by the induction hypothesis, we obtain the summands on the LHS and hence the formula for $S_{2\ell,r}$.

	\vspace{.2cm}
	Proving the formula in the theorem for even $n=2\ell$ reduces to establishing the identities
	\eqref{S2l:reven}--\eqref{S2l:rodd} below:
	\begin{align}
		\label{S2l:reven}
		[2\ell]S_{2\ell,r} &= [2\ell-r] K^{-1} S_{2\ell-1,r} +(\Y+F) S_{2\ell-1,r-1},   \text{ for } r=2m \text{ even},
	\end{align}
	and
	\begin{align}
		\label{S2l:rodd}
		& [2\ell] S_{2\ell,r} =
		[2\ell-r] K^{-1} S_{2\ell-1,r} +(\Y+F) S_{2\ell-1,r-1} +[2\ell-r+1] (\qvs)K^{-1} S_{2\ell-1,r-2},
		\\
		&\qquad\qquad\qquad\qquad\qquad \qquad\qquad\qquad\qquad\qquad \qquad\qquad\quad
		\text{ for } r=2m+1 \text{ odd}.
		\notag
	\end{align}
	
	\vspace{.3cm}
	We first prove \eqref{S2l:reven}.
	By a direct computation, we have
	\begin{align}
		\begin{split}
			q^{r-2a}[2\ell-r]  +q^{2\ell-a}[a]
			+q^{2c+r-2\ell-a} & [r-2c-a]\frac{q^{-2r}K^{-2}-1}{q^{4c-2r}K^{-2}-1}
			\label{eq:K2} \\
			&+q^{2c-2\ell}[2c] \frac{q^{-2a}K^{-2}-1}{q^{4c-2r}K^{-2}-1 }
			= [2\ell].
		\end{split}
	\end{align}
	
	Then we compute
	\begin{align*}
		\text{RHS} &\eqref{S2l:reven}\\
		=& [2\ell-r]K^{-1}\Big(\sum_{c=0}^{\frac{r}2 }\sum_{a=0}^{r-2c}  q^{\binom{2c+1}{2}+(r-2c)(r+1-2\ell)-a(r-2c-a)}(\qvs)^c \Y^{(a)}  \qbinom{\kk;1- \frac{r}2 }{c} K^{r+1-2\ell} F^{(r-2c-a)} \Big)\\
		&+(\Y+F)\Big(\sum_{c=0}^{ \frac{r}2-1 }\sum_{a=0}^{r-1-2c}  q^{\binom{2c+1}{2}+(r-1-2c)(r-2\ell-a)+a^2} (\qvs)^c\Y^{(a)}  \qbinom{\kk;1- \frac{r}2 }{c} K^{r-2\ell} F^{(r-1-2c-a)} \Big)
		\\
		=&[2\ell-r]\Big(\sum_{c=0}^{\frac{r}2 }\sum_{a=0}^{r-2c}  q^{\binom{2c+1}{2}+(r-2c)(r+1-2\ell)-a(r-2c-a)-2a} (\qvs)^c\Y^{(a)}  \qbinom{\kk;1- \frac{r}2 }{c} K^{r-2\ell} F^{(r-2c-a)} \Big)\\
		&+\sum_{c=0}^{ \frac{r}2-1 }\sum_{a=0}^{r-1-2c}  q^{\binom{2c+1}{2}+(r-1-2c)(r-2\ell-a)+a^2}[a+1](\qvs)^c\Y^{(a+1)}  \qbinom{\kk;1- \frac{r}2 }{c} K^{r-2\ell} F^{(r-1-2c-a)} \\
		&+\sum_{c=0}^{ \frac{r}2-1 }\sum_{a=0}^{r-1-2c}  q^{\binom{2c+1}{2}+(r-1-2c)(r-2\ell-a)+a^2-2a+2r-4\ell} (\qvs)^c\Y^{(a)}[r-2c-a] \qbinom{\kk;- \frac{r}2 }{c} K^{r-2\ell} F^{(r-2c-a)} \\
		&+\sum_{c=0}^{ \frac{r}2-1 }\sum_{a=0}^{r-1-2c}  q^{\binom{2c+1}{2}+(r-1-2c)(r-2\ell-a)+a^2} (\qvs)^{c+1}\Y^{(a-1)}\frac{q^{3-3a}K^{-2}-q^{1-a}}{q^2-1}  
		\\
		&\quad\quad\cdot\qbinom{\kk;1- \frac{r}2 }{c} K^{r-2\ell} F^{(r-1-2c-a)}
		\\
		=&[2\ell-r]\Big(\sum_{c=0}^{\frac{r}2 }\sum_{a=0}^{r-2c}  q^{\binom{2c+1}{2}+(r-2c)(r+1-2\ell)-a(r-2c-a)-2a}(\qvs)^c \Y^{(a)}  \qbinom{\kk;1- \frac{r}2 }{c} K^{r-2\ell} F^{(r-2c-a)} \Big)\\
		&+\sum_{c=0}^{ \frac{r}2-1 }\sum_{a=1}^{r-2c}  q^{\binom{2c+1}{2}+(r-1-2c)(r-2\ell)-(a-1)(r-2c-a)}[a](\qvs)^c\Y^{(a)}  \qbinom{\kk;1- \frac{r}2 }{c} K^{r-2\ell} F^{(r-2c-a)} \\
		&+\sum_{c=0}^{ \frac{r}2-1 }\sum_{a=0}^{r-1-2c}  q^{\binom{2c+1}{2}+(r-1-2c)(r-2\ell-a)+a^2-2a+2r-4\ell} (\qvs)^c\Y^{(a)}[r-2c-a] \qbinom{\kk;- \frac{r}2 }{c} K^{r-2\ell} F^{(r-2c-a)} \\
		&+\sum_{c=1}^{ \frac{r}2 }\sum_{a=0}^{r-2c}  q^{\binom{2c-1}{2}+(r+1-2c)(r-2\ell)-(a+1)(r-2c-a)} (\qvs)^c\Y^{(a)}\frac{q^{-3a}K^{-2}-q^{-a}}{q^2-1}  \\
		&\qquad\qquad\qquad\cdot\qbinom{\kk;1- \frac{r}2 }{c-1} K^{r-2\ell} F^{(r-2c-a)},
	\end{align*}
	which is equal to
	\begin{align*}
		=&[2\ell-r]\Big(\sum_{c=0}^{\frac{r}2 }\sum_{a=0}^{r-2c}  q^{\binom{2c+1}{2}+(r-2c)(r+1-2\ell)-a(r-2c-a)-2a}(\qvs)^c \Y^{(a)}  \qbinom{\kk;1- \frac{r}2 }{c} K^{r-2\ell} F^{(r-2c-a)} \Big)\\
		&+\sum_{c=0}^{ \frac{r}2-1 }\sum_{a=1}^{r-2c}  q^{\binom{2c+1}{2}+(r-1-2c)(r-2\ell)-(a-1)(r-2c-a)}[a](\qvs)^c\Y^{(a)}  \qbinom{\kk;1- \frac{r}2 }{c} K^{r-2\ell} F^{(r-2c-a)} \\
		& +\sum_{c=0}^{ \frac{r}2-1 }\sum_{a=0}^{r-1-2c}  q^{\binom{2c+1}{2}+(r+1-2c)(r-2\ell-a)+a^2} (\qvs)^c\Y^{(a)}[r-2c-a]\frac{q^{-2r}K^{-2}-1}{q^{4c-2r}K^{-2}-1} \qbinom{\kk;1- \frac{r}2 }{c} 
		\\
		&\cdot K^{r-2\ell} F^{(r-2c-a)}  +\sum_{c=1}^{ \frac{r}2 }\sum_{a=0}^{r-2c}  q^{\binom{2c-1}{2}+(r+1-2c)(r-2\ell)-(a+1)(r-2c-a)} (\qvs)^c\Y^{(a)}
		\\
		&\cdot\frac{(q^{-3a}K^{-2}-q^{-a})(q^{4c}-1)}{(q^2-1)(q^{4c-2r}K^{-2}-1)}\qbinom{\kk;1- \frac{r}2 }{c} K^{r-2\ell} F^{(r-2c-a)}
		\\
		&= \sum_{c=0}^{\frac{r}2 }\sum_{a=0}^{r-2c}  q^{\binom{2c}{2}+(r-2c)(r-2\ell)-a(r-2c-a)}(\qvs)^c \Y^{(a)} \\
		&\quad \cdot \left(q^{r-2a}[2\ell-r]+q^{2\ell-a}[a]+q^{2c+r-2\ell-a}[r-2c-a]\frac{q^{-2r}K^{-2}-1}{q^{4c-2r}K^{-2}-1}+q^{2c-2\ell}[2c] \frac{q^{-2a}K^{-2}-1}{q^{4c-2r}K^{-2}-1 }\right) \\
		&\quad \cdot  \qbinom{\kk;1- \frac{r}2 }{c} K^{r-2\ell} F^{(r-2c-a)}\\
		& \stackrel{\eqref{eq:K2}}{=} [2\ell]\sum_{c=0}^{\frac{r}2 }\sum_{a=0}^{r-2c}  q^{\binom{2c}{2}+(r-2c)(r-2\ell)-a(r-2c-a)} (\qvs)^c\Y^{(a)}  \qbinom{\kk;1- \frac{r}2 }{c} K^{r-2\ell} F^{(r-2c-a)} \\
		&= [2\ell] S_{2\ell,r}=\text{LHS}\eqref{S2l:reven}.
	\end{align*}

	\vspace{.3cm}
	Next we prove \eqref{S2l:rodd}. By a direct computation, we have
	\begin{align}
		\label{eq:K3}
		\begin{split}
			q^{r-2a}[2\ell-r] & \frac{q^{2-2r}K^{-2}-1}{q^{4c-2r+2}K^{-2}-1}  +q^{2\ell-a}[a]+q^{2c+r-2\ell-a}[r-2c-a]\frac{q^{2-2r}K^{-2}-1}{q^{4c-2r+2}K^{-2}-1}\\
			&+q^{2c-2\ell}[2c]\frac{q^{-2a}K^{-2}-1}{q^{4c-2r+2}K^{-2}-1}+q^{1-r-2a}[2\ell-r+1] \frac{(q^{4c}-1)K^{-2}}{q^{4c-2r+2}K^{-2}-1}
			=[2\ell].
		\end{split}
	\end{align}
	Then we  compute
	\begin{align*}
		&\text{RHS} \eqref{S2l:rodd}
		\\
		&= [2\ell-r]K^{-1}\Big(\sum_{c=0}^{\frac{r-1}2 }\sum_{a=0}^{r-2c}  q^{\binom{2c+1}{2}+(r-2c)(r+1-2\ell)-a(r-2c-a)} (\qvs)^c\Y^{(a)}  \qbinom{\kk;- \frac{r-1}2 }{c} K^{r+1-2\ell} F^{(r-2c-a)} \Big)
		\\
		& +(\Y+F)\Big(\sum_{c=0}^{ \frac{r-1}2 }\sum_{a=0}^{r-1-2c}  q^{\binom{2c+1}{2}+(r-1-2c)(r-2\ell-a)+a^2}(\qvs)^c \Y^{(a)}  \qbinom{\kk;- \frac{r-3}2 }{c} K^{r-2\ell} F^{(r-1-2c-a)} \Big)
		\\
		& +[2\ell-r+1]K^{-1}(\qvs)\Big(\sum_{c=0}^{\frac{r-3}2 }\sum_{a=0}^{r-2-2c}  q^{\binom{2c+1}{2}+(r-2-2c)(r-1-2\ell-a)+a^2} (\qvs)^c\Y^{(a)}  \\
		&\qquad\qquad\cdot\qbinom{\kk;- \frac{r-3}2 }{c} K^{r-1-2\ell} F^{(r-2-2c-a)} \Big),
	\end{align*}
	which is equal to
	\begin{align*}
		=& [2\ell-r]\Big(\sum_{c=0}^{\frac{r-1}2 }\sum_{a=0}^{r-2c}  q^{\binom{2c+1}{2}+(r-2c)(r+1-2\ell)-a(r+2-2c-a)} (\qvs)^c\Y^{(a)}\frac{q^{2-2r}K^{-2}-1}{q^{4c-2r+2}K^{-2}-1}  \qbinom{\kk;- \frac{r-3}2 }{c}\\
		&\cdot K^{r-2\ell} F^{(r-2c-a)} \Big) +\sum_{c=0}^{ \frac{r-1}2 }\sum_{a=0}^{r-1-2c}  q^{\binom{2c+1}{2}+(r-1-2c)(r-2\ell-a)+a^2}[a+1] (\qvs)^c\Y^{(a+1)}  \qbinom{\kk;- \frac{r-3}2 }{c} 
		\\
		&\cdot K^{r-2\ell} F^{(r-1-2c-a)} +\sum_{c=0}^{ \frac{r-1}2 }\sum_{a=0}^{r-1-2c}  q^{\binom{2c+1}{2}+(r+1-2c)(r-2\ell-a)+a^2}[r-2c-a] (\qvs)^c\Y^{(a)} \frac{q^{2-2r}K^{-2}-1}{q^{4c-2r+2}K^{-2}-1}   
		\\
		&\cdot\qbinom{\kk;- \frac{r-3}2 }{c} K^{r-2\ell} F^{(r-2c-a)}
		+\sum_{c=0}^{ \frac{r-3}2 }\sum_{a=0}^{r-1-2c}  q^{\binom{2c+1}{2}+(r-1-2c)(r-2\ell-a)+a^2} (\qvs)^{c+1}\Y^{(a-1)} \frac{q^{3-3a}K^{-2}-q^{1-a}}{q^2-1} \\
		&\cdot\qbinom{\kk;- \frac{r-3}2 }{c} K^{r-2\ell} F^{(r-1-2c-a)}
		+[2\ell-r+1]
		\\
		&\cdot\Big(\sum_{c=0}^{\frac{r-3}2 }\sum_{a=0}^{r-2-2c}  q^{\binom{2c+1}{2}+(r-2-2c)(r-1-2\ell)-a(r-2c-a)} (\qvs)^{c+1}\Y^{(a)}  \qbinom{\kk;- \frac{r-3}2 }{c} K^{r-2-2\ell} F^{(r-2-2c-a)} \Big)
		\\
		&= \sum_{c=0}^{\frac{r-1}2 }\sum_{a=0}^{r-2c}  q^{\binom{2c}{2}+(r-2c)(r-2\ell)-a(r-2c-a)} (\qvs)^c\Y^{(a)}\\
		&\quad \cdot\Big( q^{r-2a}[2\ell-r]\frac{q^{2-2r}K^{-2}-1}{q^{4c-2r+2}K^{-2}-1}  +q^{2\ell-a}[a]+q^{2c+r-2\ell-a}[r-2c-a]\frac{q^{2-2r}K^{-2}-1}{q^{4c-2r+2}K^{-2}-1}\\
		&\quad +q^{2c-2\ell}[2c]\frac{q^{-2a}K^{-2}-1}{q^{4c-2r+2}K^{-2}-1}+q^{1-r-2a}[2\ell-r+1] \frac{(q^{4c}-1)K^{-2}}{q^{4c-2r+2}K^{-2}-1} \Big)\qbinom{\kk;- \frac{r-3}2 }{c}
		K^{r-2\ell} F^{(r-2c-a)}
		\\
		& \stackrel{\eqref{eq:K3}}{=} [2\ell]\sum_{c=0}^{\frac{r-1}2 }\sum_{a=0}^{r-2c}  q^{\binom{2c}{2}+(r-2c)(r-2\ell)-a(r-2c-a)} (\qvs)^c\Y^{(a)}\qbinom{\kk;- \frac{r-3}2 }{c} K^{r-2\ell} F^{(r-2c-a)}\\
		& =[2\ell] S_{2\ell,r}=\text{LHS} \eqref{S2l:rodd}.
	\end{align*}
	This completes Step {\bf I}.

	\vspace{.3cm}
	{\bf Step~II.}
	Assuming the first formula in the theorem for even $n=2\ell$, we shall prove the second formula in the theorem for $n=2\ell+1$.
	It follows from \eqref{eq:ttev} that
	\begin{align}
		\label{D:22}
		[2\ell+1] \Delta(\dvev{2\ell+1}) =& \Delta(\ttt) \Delta(\dvev{2\ell})-[2\ell]\Delta(\dvev{2\ell-1})\\
		= &\left(\ttt \otimes K^{-1} + 1 \otimes (\Y +  F) \right)\cdot \Delta(\dvev{2\ell})-[2\ell]\Delta(\dvev{2\ell-1}) .\notag
	\end{align}
	We shall compare the summands $\dvev{2\ell+1-r} \otimes -$ in both sides of \eqref{D:22}. As the summands on the \text{RHS} of
	\eqref{D:22} are known by the induction hypothesis, we obtain the summands on the \text{LHS} and hence the formula for $S_{2\ell+1,r}$.

	\vspace{.2cm}
	Proving the formula in the theorem for even $n=2\ell+1$ reduces to establishing the identities
	\eqref{S22:reven}--\eqref{S22:rodd} below:
	\begin{align}
		\label{S22:reven}
		[2\ell+1]S_{2\ell+1,r} &= [2\ell+1-r] K^{-1} S_{2\ell,r}
		+[2\ell-r+2](\qvs) K^{-1} S_{2\ell,r-2} \\
		& +(\Y+F) S_{2\ell,r-1} -[2\ell](\qvs)S_{2\ell-1,r-2},  \qquad\qquad
		\text{ for } r=2m \text{ even}, \notag
	\end{align}
	and
	\begin{align}
		\label{S22:rodd}
		[2\ell+1]S_{2\ell+1,r} &= [2\ell+1-r] K^{-1} S_{2\ell,r}
		+(\Y+F) S_{2\ell,r-1}-[2\ell](\qvs)S_{2\ell-1,r-2}
		\\
		&\qquad\qquad\qquad\qquad\qquad \qquad\qquad\qquad
		\text{ for } r=2m+1 \text{ odd}.
		\notag
	\end{align}
	
	\vspace{.3cm}
	Let us first prove \eqref{S22:reven}. It follows by a direct computation that
	\begin{align}
		\label{K5}
		&q^{r-4c-2a}[2\ell+1-r]+q^{-4c-r+3-2a}\frac{(q^{4c}-1)K^{-2}}{q^{4-2r}K^{-2}-1}[2\ell+2-r]
		+q^{-4c-a+2\ell+1}\frac{q^{4c+4-2r}K^{-2}-1}{q^{4-2r}K^{-2}-1  }[a]\\
		&+q^{r-a-2c-1-2\ell}[r-2c-a]+q^{-4c+2-2\ell}\frac{q^{-2a}K^{-2}-1}{q^2-1}\frac{q^{4c}-1}{q^{4-2r}K^{-2}-1}
		-q^{-4c+1}\frac{q^{4c}-1}{q^{4-2r}K^{-2}-1}[2\ell]
		\notag \\
		&=[2\ell+1].
		\notag
	\end{align}
	Then we compute
	\begin{align*}
		&\text{LHS} \eqref{S22:reven}
		\\
		&=[2\ell+1-r] K^{-1} S_{2\ell,r}
		+[2\ell-r+2] (\qvs)K^{-1} S_{2\ell,r-2}+(\Y+F) S_{2\ell,r-1}-[2\ell](\qvs)S_{2\ell-1,r-2}\\
		&=[2\ell+1-r] K^{-1} \Big(\sum_{c=0}^{ \frac{r}2 }\sum_{a=0}^{r-2c} q^{\binom{2c}{2}+(r-2c)(r-2\ell)-a(r-2c-a)} (\qvs)^c\Y^{(a)}   \qbinom{\kk;- \frac{r-2}2 }{c} K^{r-2\ell} F^{(r-2c-a)} \Big)\\
		&\quad +[2\ell+2-r] K^{-1} (\qvs)\Big(\sum_{c=0}^{ \frac{r-2}2 }\sum_{a=0}^{r-2-2c} q^{\binom{2c}{2}+(r-2-2c)(r-2-2\ell-a)+a^2} (\qvs)^c\Y^{(a)}   \qbinom{\kk;- \frac{r-4}2 }{c}\\
		&\quad\cdot  K^{r-2-2\ell} F^{(r-2-2c-a)} \Big)+(\Y+F) \Big(\sum_{c=0}^{ \frac{r-2}2 }\sum_{a=0}^{r-1-2c} q^{\binom{2c}{2}+(r-1-2c)(r-1-2\ell-a)+a^2} (\qvs)^c\Y^{(a)}   \qbinom{\kk;- \frac{r-4}2 }{c}\\
		&\quad\cdot K^{r-1-2\ell} F^{(r-1-2c-a)} \Big)-[2\ell](\qvs)\Big(\sum_{c=0}^{ \frac{r-2}2 }\sum_{a=0}^{r-2-2c}  q^{\binom{2c+1}{2}+(r-2-2c)(r-2\ell-1-a)+a^2}(\qvs)^{c} \Y^{(a)}  
		\\
		&\qquad\qquad\cdot\qbinom{\kk;-\frac{r-4}2 }{c} K^{r-2\ell-1} F^{(r-2-2c-a)} \Big)
		\\
		& =[2\ell+1-r] \Big(\sum_{c=0}^{ \frac{r}2 }\sum_{a=0}^{r-2c} q^{\binom{2c}{2}+(r-2c)(r-2\ell)-a(r-2c-a+2)}(\qvs)^c \Y^{(a)}   \qbinom{\kk;- \frac{r-2}2 }{c}
		\\
		&\qquad\qquad\qquad\cdot K^{r-2\ell-1} F^{(r-2c-a)} \Big)+[2\ell+2-r] 
		\\
		&\qquad\cdot \Big(\sum_{c=0}^{ \frac{r-2}2 }\sum_{a=0}^{r-2-2c} q^{\binom{2c}{2}+(r-2-2c)(r-2-2\ell)-a(r-2c-a)} (\qvs)^{c+1}\Y^{(a)}   \qbinom{\kk;- \frac{r-4}2 }{c} K^{r-3-2\ell} F^{(r-2-2c-a)} \Big)
		\\
		&\quad+ \Big(\sum_{c=0}^{ \frac{r-2}2 }\sum_{a=0}^{r-1-2c} q^{\binom{2c}{2}+(r-1-2c)(r-1-2\ell-a)+a^2}[a+1] (\qvs)^c\Y^{(a+1)}   \qbinom{\kk;- \frac{r-4}2 }{c} K^{r-1-2\ell} F^{(r-1-2c-a)} \Big)
		\\
		&\quad+ \Big(\sum_{c=0}^{ \frac{r-2}2 }\sum_{a=0}^{r-1-2c} q^{\binom{2c}{2}+(r+1-2c)(r-1-2\ell-a)+a^2} [r-2c-a ](\qvs)^c\Y^{(a)}   \qbinom{\kk;- \frac{r-2}2 }{c} K^{r-1-2\ell} F^{(r-2c-a)} \Big)
		\\
		&\quad+ \Big(\sum_{c=0}^{ \frac{r-2}2 }\sum_{a=0}^{r-1-2c} q^{\binom{2c}{2}+(r-1-2c)(r-1-2\ell-a)+a^2} (\qvs)^{c+1}\Y^{(a-1)} \frac{q^{3-3a}K^{-2}-q^{1-a}}{q^2-1}  \qbinom{\kk;- \frac{r-4}2 }{c} \\
		&\qquad\qquad \cdot K^{r-1-2\ell} F^{(r-1-2c-a)} \Big)
		\\
		&\quad-[2\ell]\Big(\sum_{c=0}^{ \frac{r-2}2 }\sum_{a=0}^{r-2-2c}  q^{\binom{2c+1}{2}+(r-2-2c)(r-2\ell-1-a)+a^2}(\qvs)^{c+1} \Y^{(a)}  \qbinom{\kk;-\frac{r-4}2 }{c} K^{r-2\ell-1} F^{(r-2-2c-a)} \Big)
		\\
		& =\sum_{c=0}^{ \frac{r}2 }\sum_{a=0}^{r-2c} q^{\binom{2c+1}{2}+(r-2c)(r-2\ell-1)-a(r-2c-a)} (\qvs)^c\Y^{(a)}\cdot \Big(q^{r-4c-2a}[2\ell+1-r]
		\\
		&\quad+q^{-4c-r+3-2a}\frac{(q^{4c}-1)K^{-2}}{q^{4-2r}K^{-2}-1}[2\ell+2-r]
		+q^{-4c-a+2\ell+1}\frac{q^{4c+4-2r}K^{-2}-1}{q^{4-2r}K^{-2}-1  }[a]\\
		&\quad +q^{r-a-2c-1-2\ell}[r-2c-a]+q^{-4c+2-2\ell}\frac{q^{-2a}K^{-2}-1}{q^2-1}\frac{q^{4c}-1}{q^{4-2r}K^{-2}-1}
		\\
		&\quad-q^{-4c+1}\frac{q^{4c}-1}{q^{4-2r}K^{-2}-1}[2\ell]\Big)\cdot \qbinom{\kk;-\frac{r-2}2 }{c} K^{r-2\ell-1} F^{(r-a-2c)}
		\\
		&\stackrel{\eqref{K5}}{=} [2\ell+1]S_{2\ell+1,r}=\text{RHS}  \eqref{S22:reven}.
	\end{align*}
	
	\vspace{.3cm}
	Next we prove \eqref{S22:rodd}. By a direct computation, we have
	\begin{align}
		\label{K4}
		&q^{r-4c-2a}[2\ell+1-r]\frac{q^{4c+2-2r}K^{-2}-1}{q^{2-2r}K^{-2}-1}+
		q^{1+2\ell-4c-a}[a]\frac{q^{4c+2-2r}K^{-2}-1}{q^{2-2r}K^{-2}-1}
		\\
		&+q^{-2c+r-1-2\ell-a}[r-2c-a]
		+q^{-2c+1-2\ell}[2c]\frac{q^{-2a}K^{-2}-1}{q^{2-2r}K^{-2}-1}-q^{-4c+1}[2\ell]\frac{q^{4c}-1}{q^{2-2r}K^{-2}-1}
		\notag \\
		&=[2\ell+1].
		\notag
	\end{align}
	Then we compute
	\begin{align*}
		&\text{LHS}  \eqref{S22:rodd}
		=[2\ell+1-r] K^{-1} S_{2\ell,r}+ (\Y+F)S_{2\ell,r-1}-[2\ell](\qvs)S_{2\ell-1,r-2},
	\end{align*}
	which is equal to
	\begin{align*}
		&=[2\ell+1-r] K^{-1}    \Big(\sum_{c=0}^{ \frac{r-1}2}\sum_{a=0}^{r-2c} q^{\binom{2c}{2}+(r-2c)(r-2\ell)-a(r-2c-a)} (\qvs)^c\Y^{(a)}   \qbinom{\kk;- \frac{r-3}2}{c} K^{r-2\ell} F^{(r-2c-a)} \Big)\\
		&\quad +(\Y+F)  \Big(\sum_{c=0}^{ \frac{r-1}2 }\sum_{a=0}^{r-1-2c} q^{\binom{2c}{2}+(r-1-2c-a)(r-1-2\ell)+a^2} (\qvs)^c\Y^{(a)}   \qbinom{\kk;-\frac{r-3}2 }{c} K^{r-1-2\ell} F^{(r-1-2c-a)} \Big)\\
		&\quad -[2\ell] \Big(\sum_{c=0}^{ \frac{r-3}2 }\sum_{a=0}^{r-2-2c}  q^{\binom{2c+1}{2}+(r-2-2c)(r-1-2\ell-a)+a^2}(\qvs)^{c+1} \Y^{(a)}  \qbinom{\kk;- \frac{r-3}2 }{c} K^{r-1-2\ell} F^{(r-2-2c-a)} \Big)
		\\
		&=[2\ell+1-r]     \Big(\sum_{c=0}^{ \frac{r-1}2}\sum_{a=0}^{r-2c} q^{\binom{2c}{2}+(r-2c)(r-2\ell)-a(r-2c-a+2)} (\qvs)^c\Y^{(a)}   \qbinom{\kk;- \frac{r-3}2}{c} K^{r-2\ell-1} F^{(r-2c-a)} \Big)
		\\
		&\quad+\Big(\sum_{c=0}^{ \frac{r-1}2 }\sum_{a=0}^{r-1-2c} q^{\binom{2c}{2}+(r-1-2c)(r-1-2\ell-a)+a^2}[a+1] (\qvs)^c\Y^{(a+1)}   \qbinom{\kk;-\frac{r-3}2 }{c} K^{r-1-2\ell} F^{(r-1-2c-a)} \Big)
		\\
		&\quad + \Big(\sum_{c=0}^{ \frac{r-1}2 }\sum_{a=0}^{r-1-2c} q^{\binom{2c}{2}+(r-1-2c)(r-1-2\ell-a)+a^2+2a}(\qvs)^c \Y^{(a)}   \qbinom{\kk;-\frac{r-1}2 }{c} q^{2r-2-4\ell}[r-2c-a]
		\\
		&\qquad\qquad \cdot K^{r-1-2\ell} F^{(r-2c-a)} \Big)
		\\
		&\quad + \Big(\sum_{c=0}^{ \frac{r-1}2 }\sum_{a=0}^{r-1-2c} q^{\binom{2c}{2}+(r-1-2c)(r-1-2\ell-a)+a^2}(\qvs)^{c+1} \Y^{(a-1)}  \frac{q^{3-3a}K^{-2}-q^{1-a}}{q^2-1}\qbinom{\kk;-\frac{r-3}2 }{c}
		\\
		&\qquad\qquad \cdot K^{r-1-2\ell} F^{(r-1-2c-a)} \Big)
		\\
		&\quad -[2\ell] \Big(\sum_{c=0}^{ \frac{r-3}2 }\sum_{a=0}^{r-2-2c}  q^{\binom{2c+1}{2}+(r-2-2c)(r-1-2\ell-a)+a^2}(\qvs)^{c+1} \Y^{(a)}  \qbinom{\kk;- \frac{r-3}2 }{c} K^{r-1-2\ell} F^{(r-2-2c-a)} \Big)
		\\
		&=\sum_{c=0}^{ \frac{r-1}2}\sum_{a=0}^{r-2c} q^{\binom{2c+1}{2}+(r-2c)(r-2\ell-1-a)+a^2} (\qvs)^c\Y^{(a)}\cdot\Big(q^{r-4c-2a}[2\ell+1-r]\frac{q^{4c+2-2r}K^{-2}-1}{q^{2-2r}K^{-2}-1}  \\
		&+ q^{1+2\ell-4c-a}[a]\frac{q^{4c+2-2r}K^{-2}-1}{q^{2-2r}K^{-2}-1} +q^{-2c+r-1-2\ell-a}[r-2c-a]\\
		&+q^{-2c+1-2\ell}[2c]\frac{q^{-2a}K^{-2}-1}{q^{2-2r}K^{-2}-1}-q^{-4c+1}[2\ell]\frac{q^{4c}-1}{q^{2-2r}K^{-2}-1}   \Big)
		\qbinom{\kk;-\frac{r-1}2 }{c} K^{r-1-2\ell} F^{(r-2c-a)}\\
		&\stackrel{\eqref{K4}}{=} [2\ell+1]S_{2\ell+1,r}=\text{RHS}  \eqref{S22:rodd}.
	\end{align*}
	
	This completes the proof of Theorem \ref{thm:even}.

	\subsection{A new form of the formula for $\Delta(\dvev{n})$}
	
	Let us reformulate the formulae in Theorem \ref{thm:even} in another form.
	
	\begin{prop}
		\label{prop:FHY}
		For $n$ even, we have
		{\small
			\begin{align*}  
				\Delta(\dvev{n})  &=
				\sum_{r=0}^n \dvev{n-r}  \otimes
				\Big(\sum_{c=0}^{\lfloor \frac{r}2 \rfloor}\sum_{a=0}^{r-2c} (-1)^cq^{3c-(r-2c)(r-n)+a(r-2c-a)}(\qvs)^c F^{(a)}   \qbinom{\kk;1-c+\lfloor \frac{r-2}2 \rfloor}{c} K^{r-n} \Y^{(r-2c-a)} \Big).
			\end{align*}
		}
		For $n$ odd, we have
		{\small
			\begin{align*}
				\Delta(\dvev{n})  &=
				\sum_{r=0}^n \dvev{n-r}  \otimes
				\Big(\sum_{c=0}^{\lfloor \frac{r}2 \rfloor}\sum_{a=0}^{r-2c} (-1)^c q^{c-(r-2c)(r-n)+a(r-2c-a)} (\qvs)^c F^{(a)}  \qbinom{\kk;1-c+\lfloor \frac{r-1}2 \rfloor}{c} K^{r-n} \Y^{(r-2c-a)} \Big).
			\end{align*}
		}
	\end{prop}
	
	\begin{proof}
		Let us derive the first formula only from Theorem \ref{thm:even}, and skip the entirely similar proof of the second one. 
		
		First, we assume $\vs=q^{-1}$. Recall the anti-involution $\chi$ such that $\chi(\dvev{n})= \dvev{n}$ from Lemma~ \ref{lem:antiinvolution}. By checking on generators via the comultiplication formula \eqref{eq:Delta}, we see that $(\chi \otimes \chi) \circ \Delta \circ \chi =\Delta$. Therefore, $(\chi \otimes \chi) \Delta (\dvev{n}) = \Delta (\dvev{n})$. We compute
		\begin{align*}
			(\chi &\otimes \chi) (\text{RHS}(\ref{eq:deltaevev}) )\\
			=&
			\sum_{r=0}^n \dvev{n-r}  \otimes  \chi \Big(\sum_{c=0}^{\lfloor \frac{r}2 \rfloor}\sum_{a=0}^{r-2c} q^{\binom{2c}{2}+(r-2c)(r-n)-a(r-2c-a)} \Y^{(a)}   \qbinom{\kk;-\lfloor \frac{r-2}2 \rfloor}{c} K^{r-n} F^{(r-2c-a)} \Big) \\
			=&\sum_{r=0}^n \dvev{n-r}  \otimes   \Big(\sum_{c=0}^{\lfloor \frac{r}2 \rfloor}\sum_{a=0}^{r-2c} q^{-\binom{2c}{2}-(r-2c)(r-n-a)-a^2}
			F^{(r-2c-a)}(-1)^cq^{2c(c+1)}\qbinom{\kk;1+\lfloor \frac{r-2}2 \rfloor-c}{c}K^{r-n}
			\Y^{(a)}    \Big)\\
			=& \sum_{r=0}^n \dvev{n-r}  \otimes   \Big(\sum_{c=0}^{\lfloor \frac{r}2 \rfloor}\sum_{a=0}^{r-2c} (-1)^cq^{3c-(r-2c)(r-n-a)-a^2}
			F^{(r-2c-a)}\qbinom{\kk;1+\lfloor \frac{r-2}2 \rfloor-c}{c}K^{r-n}
			\Y^{(a)}    \Big)\\
			=& \sum_{r=0}^n \dvev{n-r}  \otimes   \Big(\sum_{c=0}^{\lfloor \frac{r}2 \rfloor}\sum_{a=0}^{r-2c} (-1)^cq^{3c-(r-2c)(r-n-a)-a^2}
			F^{(a)}\qbinom{\kk;1+\lfloor \frac{r-2}2 \rfloor-c}{c}K^{r-n}
			\Y^{(r-2c-a)}    \Big).
		\end{align*}
	This proves the first formula for $\vs=q^{-1}$, i.e.,
{\small
			\begin{align}  
				&\Delta(\dvev{n})  =
				\sum_{r=0}^n \dvev{n-r}  \otimes
				\notag \\&
				\Big(\sum_{c=0}^{\lfloor \frac{r}2 \rfloor}\sum_{a=0}^{r-2c} (-1)^cq^{3c-(r-2c)(r-n)+a(r-2c-a)} F^{(a)}   \qbinom{\kk;1-c+\lfloor \frac{r-2}2 \rfloor}{c} K^{r-n} (q^{-1}EK^{-1})^{(r-2c-a)} \Big).
				\label{special}
			\end{align}
		}
		
	We shall reduce the first formula for general parameter $\vs$ to the identity \eqref{special}. 
	We continue to denote by $B=F +q^{-1} EK^{-1}$ and $\dvev{n}$ the corresponding $\imath$divided powers with {\em special parameter} $q^{-1}$. Below we shall denote by $'B=F + \vs EK^{-1}$ and $'\dvev{n}$ the corresponding $\imath$divided powers with {\em general parameter} $\vs$. 
	
	 Below we work with $\U$ over an extension field of $\mathbb Q(q)$.  Consider the rescaling automorphism of Hopf algebra 
	\[
	\Phi_z: \U\rightarrow \U \otimes \U,
	\qquad
	F \mapsto z^{-1} F, \;\;
	E \mapsto z E, \;\;
	K \mapsto K, \;\;
	\text{ where } z:= (\qvs)^{\frac12}.
	\]
Note that 
	\begin{align}
	\label{z}
	\Phi_z(F) =z^{-1} F, 
	\quad
	\Phi_z(q^{-1}EK^{-1}) =z^{-1} \Y, 
	\quad
	\Phi_z (B) =z^{-1}\, {}'B,
	\quad
	\Phi_z (\dvev{n}) =z^{-n} \, {}' \dvev{n}, \forall n. 
	\end{align}
	Indeed, the second identity follows by $\Phi_z(q^{-1}EK^{-1}) = z q^{-1} EK^{-1} =z^{-1} (\vs EK^{-1}) =z^{-1} \Y$, and the remaining identities can be proved similarly.
	
	Note that $(\Phi_z \otimes \Phi_z) \circ \Delta =\Delta \circ \Phi_z$. Now applying $\Phi_z \otimes \Phi_z$ to both sides of \eqref{special} and using  \eqref{z}, we have obtained the first formula for general parameter $\vs$ in the proposition. 
	\end{proof}
	
	\section{Comultiplication formulae for $\dv{n}$}
	
	In this section we present the counterparts of Theorem~\ref{thm:even} and Proposition~\ref{prop:FHY} for the remaining family of $\imath$-divided powers $\dv{n}$. We shall skip the entirely similar proofs.
	
	\begin{thm}
		\label{thm:odd}
		For $n$ even, we have
		\begin{align*}  
			\Delta(\dv{n})  &=
			\sum_{r=0}^n \dv{n-r}  \otimes   \Big(\sum_{c=0}^{\lfloor \frac{r}2 \rfloor}\sum_{a=0}^{r-2c} q^{\binom{2c+1}{2}+(r-2c)(r-n)-a(r-2c-a)} (\qvs)^c \Y^{(a)}  \qbinom{\kk;-\lfloor \frac{r-1}2 \rfloor}{c} K^{r-n} F^{(r-2c-a)} \Big).
		\end{align*}
		For $n$ odd, we have
		\begin{align*}
			\Delta(\dv{n})  &=
			\sum_{r=0}^n \dv{n-r}  \otimes   \Big(\sum_{c=0}^{\lfloor \frac{r}2 \rfloor}\sum_{a=0}^{r-2c} q^{\binom{2c}{2}+(r-2c)(r-n)-a(r-2c-a)} (\qvs)^c \Y^{(a)}  \qbinom{\kk;-\lfloor \frac{r-2}2 \rfloor}{c} K^{r-n} F^{(r-2c-a)} \Big).
		\end{align*}
	\end{thm}
	\begin{example}
		\begin{align*}
			\Delta(\dv{2}) &=
			\dv{2} \otimes K^{-2}  + \ttt \otimes (q^{-1} \Y K^{-1} +q^{-1} K^{-1}F)
			+1\otimes (\Y^{(2)}  + q^{-1} \Y F+ F^{(2)}  +  q^{3} (\qvs) [\kk;0]).
			\\
			\Delta(\dv{3}) &=  \dv{3} \otimes K^{-3}  +\dv{2} \otimes (q^{-2} \Y K^{-2} +q^{-2} K^{-2}F)
			\\
			&\quad +\ttt \otimes \left(q^{-2} \Y^{(2)} K^{-1} + q^{-3} \Y K^{-1} F+ q^{-2} K^{-1} F^{(2)}  +  q (\qvs)[\kk;0] K^{-1} \right)
			\\
			&\quad +1 \otimes  \left(\Y^{(3)} + q^{-2} \Y^{(2)}  F+ q^{-2} \Y F^{(2)}  +F^{(3)} +q(\qvs) \Y [\kk;0] +q(\qvs) [\kk;0]F \right).
		\end{align*}
	\end{example}
	
	\begin{prop}
		For $n$ even, we have
		{\small
			\begin{align*}
				\Delta(\dv{n})  &=
				\sum_{r=0}^n \dv{n-r}  \otimes
				\bigg(\sum_{c=0}^{\lfloor \frac{r}2 \rfloor}\sum_{a=0}^{r-2c} q^{c-(r-2c)(r-n)+a(r-2c-a)} (\qvs)^cF^{(a)}  \qbinom{\kk;1-c+\lfloor \frac{r-1}2 \rfloor}{c} K^{r-n} \Y^{(r-2c-a)} \bigg).
			\end{align*}
		}
		For $n$ odd, we have
		{\small
			\begin{align*}
				\Delta(\dv{n})  &=
				\sum_{r=0}^n \dv{n-r}  \otimes
				\bigg(\sum_{c=0}^{\lfloor \frac{r}2 \rfloor}\sum_{a=0}^{r-2c} q^{3c-(r-2c)(r-n)+a(r-2c-a)} (\qvs)^cF^{(a)}  \qbinom{\kk;1-c+\lfloor \frac{r-2}2 \rfloor}{c} K^{r-n} \Y^{(r-2c-a)} \bigg).
			\end{align*}
		}
	\end{prop}

	
\end{document}